\documentclass[12pt,reqno,a4paper]{amsart}
\usepackage{extsizes}
\usepackage{blindtext}
\usepackage{fullpage}
\usepackage{mathtools}
\usepackage{longtable}
\usepackage{amsmath,amssymb,amsthm}
\usepackage{amscd}
\usepackage{bm}
\usepackage{hyperref}
\usepackage{dsfont}
\usepackage{enumerate}
\usepackage{epsfig}
\usepackage{float, graphicx}
\usepackage{latexsym, amsxtra}
\usepackage{mathrsfs}
\usepackage{multicol}
\usepackage[normalem]{ulem}
\usepackage{psfrag}
\usepackage[parfill]{parskip}
\usepackage{stmaryrd}
\usepackage{tikz}
\usepackage[T1]{fontenc}
\usepackage{url}
\usepackage{verbatim}
\usepackage{indentfirst}
\usepackage{tikz-cd}

\usepackage{svg}

\usepackage[all,cmtip]{xy}
\usepackage{graphicx}

\usepackage{mathtools}

\DeclareFontFamily{U}{dutchcal}{\skewchar\font=45 }
\DeclareFontShape{U}{dutchcal}{m}{n}{<-> s*[1.0] dutchcal-r}{}
\DeclareFontShape{U}{dutchcal}{b}{n}{<-> s*[1.0] dutchcal-b}{}

\DeclareMathAlphabet{\mathlcal}{U}{dutchcal}{m}{n}
\SetMathAlphabet{\mathlcal}{bold}{U}{dutchcal}{b}{n}

\makeatletter
\def\thm@space@setup{%
	\thm@preskip=2ex \thm@postskip=2ex
}
\makeatother

\newcommand{\Id}{{\mathrm{Id}}}
\hypersetup{hidelinks}

\newtheorem{thm}{Theorem~}[section]
\newtheorem{lem}[thm]{Lemma~}

\newtheorem{prop}[thm]{Proposition~}
\newtheorem{ques}[thm]{Question~}
\newtheorem{cor}[thm]{Corollary~}
\newtheorem{assumption}[thm]{Assumption~}

\theoremstyle{remark}
\newtheorem{rmk}[thm]{Remark~}
\newtheorem{ex}[thm]{Example~}

\theoremstyle{definition}
\newtheorem{defn}[thm]{Definition~}

\newcommand{\CC}{\mathbb{C}}
\newcommand{\ZZ}{\mathbb{Z}}

\newcommand{\QQ}{\mathbb{Q}}

\newcommand{\calO}{\mathcal{O}}

\newcommand{\calL}{\mathcal{L}}
\newcommand{\calZ}{\mathcal{Z}}
\newcommand{\calX}{\mathcal{X}}

\newcommand\Ima{\mathrm{Im}}

\newcommand\id{\mathrm{id}}

\newcommand{\rank}{\mathrm{rank}}

\newcommand{\Gr}{\mathrm{Gr}}

\newcommand{\Sing}{\mathrm{Sing}}
\newcommand{\Hlim}{\mathrm{H}_{\mathrm{lim}}}

\begin{document}
\title{Rigidity Criterion for Certain Calabi--Yau families}
\vspace{1.2cm}

\author{Ruiran Sun}
 \address{School of Mathematical Sciences, Xiamen University, Xiamen 361005, China}
\email{ruiransun@xmu.edu.cn}

\author{Chenglong Yu}
 \address{Center for Mathematics and Interdisciplinary Sciences, Fudan University and
Shanghai Institute for Mathematics and Interdisciplinary Sciences (SIMIS), Shanghai, China}
\email{yuchenglong@simis.cn}

\author[Kang Zuo]{Kang Zuo}
\address{School of Mathematics and Statistics, Wuhan University, Luojiashan, Wuchang, Wuhan, Hubei, 430072, P.R. China; Institut f\"ur Mathematik, Universit\"at Mainz, Mainz, Germany, 55099}
\email{zuok@uni-mainz.de}
\begin{abstract}
We prove a new rigidity criterion for families of polarized Calabi--Yau manifolds. If, near a boundary point, the total space is smooth, the relative canonical bundle is trivial, and the boundary fiber has only isolated singularities whose summed mixed Hodge spectrum is concentrated, then the family is rigid. This includes degenerations with only ordinary double points and cusps. The proof combines a local vanishing-cycle analysis with a global tensor-product decomposition of the associated variation of Hodge structures. We also prove an infinitesimal version. An explicit non-isotrivial, nonrigid family shows that the isolated-singularity assumption alone is insufficient.
\end{abstract}

\subjclass[2010]{14D22,14C30}
\keywords{}

\maketitle

\setcounter{tocdepth}{1}
\tableofcontents

\section{Introduction}
In this paper, we investigate the rigidity of families of Calabi--Yau manifolds over complex algebraic curves. The study of rigidity in algebraic geometry can be traced back to the classical Shafarevich conjecture. While rigidity is well-established for families of hyperbolic curves, the situation for higher-dimensional varieties, particularly Calabi--Yau manifolds, is significantly more subtle. Our primary goal is to provide a criterion for rigidity based on the nature of the singularities developed by the family at the boundary of its compactification.

\subsection*{A rigidity criterion}
Let $f\colon (\calX, \calL)\to S$ be a flat family of $n$-dimensional polarized smooth projective varieties over a quasiprojective base $S$. We say the family is \textbf{nonrigid} if there exists a quasiprojective variety $B$ with a base point $b$ and a flat polarized family $\widetilde{f}\colon (\widetilde{\calX}, \widetilde{\calL})\to S\times B$ such that the restriction $(\widetilde{\calX}, \widetilde{\calL})|_{S\times \{b\}} \cong (\calX, \calL)$ as an $S$-scheme, and the total family $(\widetilde{\calX}, \widetilde{\calL})$ is non-isotrivial in the $B$ direction. Conversely, if no such deformation exists, we say the family $f$ is \textbf{rigid}.

The motivation for this work stems from the observation that known examples of nonrigid Calabi--Yau families, such as those constructed by Viehweg and Zuo \cite{viehweg2005complex, viehweg2005geometry}, tend to develop singular fibers with singular loci of positive dimension. Isolated singularities alone, however, do not imply rigidity, even for a non-isotrivial family; see Example~\ref{example:nonrigid-cone} below. Our criterion requires, in addition, a concentration condition on the sum of the local mixed Hodge spectra.

\begin{assumption}
\label{assum: main}
Let $f\colon (\calX, \calL)\to S$ be a flat family of $n$-dimensional polarized smooth Calabi--Yau manifolds over a quasiprojective curve $S$. Let $\overline{S}$ be a projective smooth compactification of $S$ and let $0 \in \overline{S} \setminus S$ be a boundary point. We assume that there is a sufficiently small analytic disk $\Delta\subset \overline S$ centered at $0$ such that the restriction of the family to $\Delta^*=\Delta\setminus\{0\}$ extends to a flat projective morphism
\[
\overline f_\Delta:\overline{\calX}_{\Delta}\longrightarrow \Delta
\]
with the following properties:
\begin{enumerate}
\item the total space $\overline{\calX}_{\Delta}$ is smooth,
\item the relative canonical bundle $\omega_{\overline{\calX}_{\Delta}/\Delta}\cong \calO_{\overline{\calX}_{\Delta}}$ is trivial,
\item the central fiber $X_0=\overline f_\Delta^{-1}(0)$ is singular and all its singularities are isolated, say $\Sing(X_0)=\{p_1,\ldots,p_k\}$.
\end{enumerate}
Since the total space is smooth, these are isolated hypersurface singularities.
\end{assumption}

Our main result is the following rigidity criterion.

\begin{thm}\label{main-thm}
In the setting of Assumption \ref{assum: main}, suppose that the summed mixed Hodge spectrum
\[
\sigma(X_0)=\sum_{i=1}^k\sigma(X_0,p_i)
\]
is \textbf{concentrated}. Then the family $f \colon \calX \to S$ is rigid.
\end{thm}
See Definition \ref{defn: concentrated mixed Hodge spectrum} for the precise definition of concentrated mixed Hodge spectrum.
As a direct application, since any collection of ordinary double points (ODP) and cusp singularities in a fixed dimension has concentrated summed spectrum (Example~\ref{example: A-D-E}), we obtain:

\begin{cor}\label{main-cor}
In the setting of Assumption \ref{assum: main}, if every singularity of $X_0$ is an ordinary double point or a cusp singularity, then the family $f \colon \calX \to S$ is rigid.
\end{cor}

The same spectral condition also forces the $(-1,1)$-part of the flat endomorphism algebra of the Calabi--Yau isotypic component to vanish. We prove this infinitesimal statement in Theorem~\ref{thm:infinitesimal} and recover the rigidity of $f$ in Corollary~\ref{cor:infinitesimal-rigidity}.

\subsection*{Previous works}
As aforementioned, for the case of non-isotrivial hyperbolic curves, the family is always rigid. This is part of the geometric version of Shafarevich's conjecture, proved by Parshin \cite{parshin1968} and Arakelov \cite{arakelov1971}.

In higher dimensions, the behavior is more complicated due to the existence of nonrigid families. Faltings initially established a rigidity criterion for families of abelian varieties \cite{faltings83}, while Peters investigated the problem for abstract variations of Hodge structures (VHS) \cite{peters90}. Building on these results, Saito-Zucker \cite{saito-zucker} and Saito \cite{saito-ab} classified non-rigid families of K3 surfaces and abelian varieties by analyzing the endomorphism algebras of the corresponding VHS.

These classification results motivated further study of the distribution of non-rigid families in moduli spaces satisfying the Torelli property \cite{CHSZ}. A corollary is that a non-isotrivial family is rigid if its classifying map covers a very general point in the moduli space. The main theorem of the present paper provides a refinement of this ``very general rigidity'' specifically for the Calabi--Yau case.

In the Calabi--Yau setting, Liu-Todorov-Yau-Zuo \cite{LTYZ05} established a rigidity criterion based on the non-vanishing of the Yukawa coupling. Subsequently, Zhang \cite{zhang04} proved that a family is rigid if it has a boundary point with \textbf{maximum unipotent monodromy} (i.e. large complex structure limit).

On the other hand, non-rigid Calabi--Yau families do exist. Viehweg-Zuo \cite{viehweg2005geometry} constructed a nonrigid family of quintic threefolds, and additional examples are provided in our forthcoming work \cite{sun2026nonrigidCY}. Notably, in these non-rigid examples, the VHS $R^n\widetilde{f}_*(\CC)$ on $S\times B$ usually admits a decomposition arising from a {\it geometric correspondence} between the $n$-dimensional family and a product of families of lower dimensional varieties.

\subsection*{A conjectural ``motivic decomposition''}

The known examples of non-rigid families inspire us to ask the following question:
\begin{ques}[Motivic decomposition]
  \label{conj: mot_decomp}
  Consider a flat family of $n$-dimensional polarized manifolds over a product surface, $f\colon (\calX, \calL)\to S \times B$, whose fibers satisfy the local Torelli property. Here $S$ and $B$ are smooth quasi-projective curves with projective completions $\overline{S}$ and $\overline{B}$. Suppose the family is non-isotrivial along both $S$ and $B$ directions. Assume further that the local system $\mathbb{V}:=R^nf_*\mathbb{C}_{\calX}$ has nontrivial monodromy along the boundary divisor $(S \times (\overline{B} \setminus B)) \cup ((\overline{S} \setminus S) \times B)$.

  Then by Deligne's theorem, each irreducible factor of $\mathbb{V}$ decompose as tensor product of VHS pulled back from $S$ and $B$ respectively. Are these tensor-product decompositions {\rm motivic}? More precisely, can one find some product family $(g,h)\colon \mathcal{Y}_1 \times \mathcal{Y}_2 \to S' \times B'$ ($S'$ and $B'$ are finite coverings of $S$ and $B$) and a correspondence between the base change of $f\colon \calX \to S \times B$ and the product family $(g,h)$ such that these tensor-products can be realized as irreducible factors of $R^{\cdot}g_*\mathbb{C}_{\mathcal{Y}_1} \otimes R^{\cdot}h_*\mathbb{C}_{\mathcal{Y}_2}$?
\end{ques}
Recall that in Faltings' examples of non-rigid families of \textit{simple} abelian varieties, the base curve is always a \textit{compact} Shimura curve; consequently, the family has no singular fiber. This explains why the condition requiring non-trivial monodromy along the boundary divisor is necessary in the statement of Question~\ref{conj: mot_decomp}.

The tensor-product decomposition on VHS imposes restriction on the Mumford-Tate group of $\mathbb{V}$. This gives some Hodge cycle $C$, which, by assuming the Hodge conjecture, corresponds to an algebraic cycle in some self-fiber product of the family $f$. However, it is not clear how $C$ might induce a geometric decomposition of the family.

In another cotext, Viehweg-Zuo use the maximal Higgs field plus the nontrivial local monodromy action along the boundary to decompose the family of abelian varieties into self-product of universal elliptic curves (up to isogeny, cf. {\cite[Theorem~0.2]{VZ04Ara}}). Motivated by these precedents, we intend to investigate Question~\ref{conj: mot_decomp} for families with high-weight VHS in the forthcoming work \cite{sun2026nonrigidCY}.

Under the assumption of Question~\ref{conj: mot_decomp}, one factor of the predicted \textit{product structure} in a non-rigid family provides the ``room'' for the deformation to occur, while the other factor remains fixed over the base $S$. Our criterion shows that, when all boundary singularities are isolated, concentration of their summed mixed spectrum is incompatible with the multiple Hodge levels required by such a non-trivial tensor product decomposition. Example~\ref{example:nonrigid-cone} explains why a condition on the singularities beyond isolation is needed.

\subsection*{Methods}
The proof proceeds via a contradiction. Suppose the family is nonrigid. By a result of Deligne \cite{deligne1987theoreme}, the irreducible factor of the VHS containing the $(n,0)$-part must decompose as a tensor product $\mathbb{V}_0\cong p_1^*\mathbb{V}_1\otimes p_2^*\mathbb{V}_2$. This decomposition persists for the limiting mixed Hodge structures (MHS) at the boundary point.

On the other hand, we analyze the local behavior near the isolated singularities using vanishing cycles. A key step is Proposition~\ref{non-triv-T}: the restriction of the holomorphic top form $\Omega_t$ to each Milnor fiber represents a nonzero cohomology class for $0<|t|\ll1$. The proof uses the residue description of relative top forms and the freeness of the Brieskorn lattice. For the nodal case, we also discuss the period-integral approach related to \cite{rollenske2009smoothing,HeinSun2017,collins2023special}. Since all singularities are isolated, the kernel of the total vanishing-cycle map is exactly the monodromy-invariant subspace. Its image on the Calabi--Yau factor therefore retains the whole second tensor factor, contradicting concentration of the summed mixed spectrum. The same kernel identity gives a faithful action of the relevant flat endomorphism algebra on the vanishing image, from which the infinitesimal version follows.

\subsection*{Organization of the Paper}
The paper is organized as follows. In Section~\ref{sec-2}, we recall the vanishing-cycle sequence and define the concentrated mixed Hodge spectrum for isolated singularities. In Section~\ref{sec-3}, we prove Proposition~\ref{non-triv-T} using a local residue argument and then discuss the period-integral approach in the ordinary-double-point case. Section~\ref{sec-4} gives the tensor-product decomposition of the limiting MHS, proves Theorem~\ref{main-thm}, and ends with a nonrigid example. In Section~\ref{sec-5}, we prove the infinitesimal version and deduce rigidity from it. Finally, in Section~\ref{sec-6}, we discuss the relation of our results to C.-L. Wang's work on the incompleteness of the Weil--Petersson metric.

\subsection*{Notations and Conventions}
Throughout this paper, we work over the field of complex numbers $\CC$. A Calabi--Yau manifold is defined as a smooth projective variety with trivial canonical bundle $K_X \cong \mathcal{O}_X$. For a family $f: \calX \to \Delta$, $\Hlim^n(X_t)$ denotes the cohomology of a nearby fiber endowed with the limiting mixed Hodge structure. We denote the semisimple and nilpotent parts of the monodromy $T$ by $T^{ss}$ and $N$, respectively.

\subsection*{Acknowledgment}
The first and second authors are grateful to Dingxin Zhang for helpful discussions on the Hodge theory of isolated singularities. The second author is supported by the national key research and development program of China (No. 2022YFA1007100) and NSFC 12201337. The third author is supported by the Key Program (Grant No. 12331002) and the International Collaboration Fund (Grant No. W2441003) of the National Natural Science Foundation of China.

\section{Vanishing cycle sequences and singularities}\label{sec-2}
We first study the local behaviour of the family near a singular fiber. In this section, let $f\colon \calX\to \Delta$ be a flat family of $n$-dimensional polarized projective varieties over the unit disk $\Delta=\{t\in \CC\mid |t|<1\}$, smooth over $\Delta-\{0\}$. Let $p\in X_0=f^{-1}(0)$ be an isolated singular point. Choosing a Milnor ball $B_p$ around $p$, we write
\[
F_p=X_t\cap B_p,\qquad 0<|t|\ll 1,
\]
for the local Milnor fiber. Restriction of nearby cycles to $B_p$ gives a local vanishing-cycle morphism of mixed Hodge structures
\begin{equation}
	\label{equation: local vanishing map}
	\nu_p:\Hlim^n(X_t)\longrightarrow \mathrm{H}^n(F_p).
\end{equation}
For the rigidity argument we use the sum of these local maps over all singularities.

If, in addition, all singularities of $X_0$ are isolated, say $p_1,\ldots,p_k$, put $H_\Sigma=\bigoplus_{i=1}^k\mathrm{H}^n(F_{p_i})$ and denote the total map by $\nu_\Sigma=\bigoplus_i\nu_{p_i}$. It fits into the usual vanishing-cycle exact sequence, see \cite{KerrLaza2021},
\begin{equation}
	\label{equation: vanishing sequences}
	0\to \mathrm{H}^n(X_0)\xrightarrow{sp} \Hlim^n(X_t)\xrightarrow{\nu_\Sigma} \bigoplus_{i=1}^k \mathrm{H}^{n}(F_{p_i})\to \mathrm{H}^{n+1}(X_0)\to \Hlim^{n+1}(X_t)\to 0.
\end{equation}
Let $T$ be the monodromy operator arising from $\pi_1(\Delta-\{0\})$. Then sequence \eqref{equation: vanishing sequences} is $T$-equivariant and respects the mixed Hodge structures. If the total space $\calX$ is smooth and all singularities of $X_0$ are isolated, then Kerr--Laza \cite[Corollary 5.7]{KerrLaza2021} prove that the morphism $sp$ induces an isomorphism
\begin{equation*}
sp\colon \mathrm{H}^n(X_0)\to \Hlim^n(X_t)^T,
\end{equation*}
where $\Hlim^n(X_t)^T$ is the $T$-invariant subspace. Consequently,
\begin{equation}\label{eq:total-kernel}
\ker\nu_\Sigma=\Hlim^n(X_t)^T.
\end{equation}
Let $T=T^{ss}T^{un}$ be the Jordan decomposition, and put $N=\log T^{un}$. Thus $N$ is nilpotent and has Hodge type $(-1,-1)$ on the limiting MHS.

In the following, we specialize to the Calabi--Yau situation in Assumption~\ref{assum: main}: locally over the boundary disk the total space is smooth and $\omega_{\calX/\Delta}$ is generated by a relative volume form. The summed mixed Hodge spectrum of all isolated singularities will provide the obstruction to deformations.

\subsection{Mixed spectrum for isolated singularities}
First we need the following definition of mixed Hodge spectrum of isolated singularities, see \cite[Definition 1.1]{KerrLazaToAppear}.
\begin{defn}
	\label{defn: Hodge spectrum}
Let $(X, x)$ be a germ of isolated singularity of dimension $n$. Let $F_x$ be the Milnor fiber of the singularity and $T$ be the monodromy operator on $\mathrm{H}^n(F_x)$. Let $T=T^{ss}T^{un}$ be the Jordan decomposition of $T$. Then $\mathrm{H}^n(F_x)$ has eigendecomposition under $T^{ss}$ as $\mathrm{H}^n(F_x) = \bigoplus_{\lambda} \mathrm{H}^n_\lambda(F_x)$ with eigenvalues $e^{2\pi \sqrt{-1} \lambda}$ for $\lambda\in [0,1)$. The mixed Hodge structure on $\mathrm{H}^n(F_x)$ gives rise to the Hodge--Deligne decomposition $\mathrm{H}^n(F_x)_{\lambda}=\bigoplus_{p,q}\mathrm{H}^{p,q}_\lambda(F_x)$. The mixed Hodge spectrum of the singularity $(X,x)$ is defined as formal spectrum
\[
\sigma(X,x)=\sum m_{\alpha, w}[(\alpha, w)]\in \ZZ[\QQ\times \ZZ]
\]
where $m_{\alpha, w} = \dim \mathrm{H}^{p,q}_\lambda(F_x)$ with $\alpha=\lambda+p$ and $w=p+q$.

If $X_0$ has several isolated singular points $p_i$, one may also consider the sum
\[
\sigma(X_0) = \sum_i \sigma(X_0, p_i).
\]
The main theorem uses this sum over all singular points of $X_0$.
\end{defn}

\begin{defn}
	\label{defn: concentrated mixed Hodge spectrum}
We say that an isolated singularity $(X,x)$ has concentrated mixed Hodge spectrum if, for each pair $(\lambda,w)\in [0,1)\times \ZZ$, there is at most one integer $p$ such that
\[
m_{\lambda+p,w}\neq 0
\]
in the mixed Hodge spectrum $\sigma(X,x)$. Equivalently, for fixed monodromy eigenvalue $e^{2\pi\sqrt{-1}\lambda}$ and fixed weight $w$, the Hodge--Deligne decomposition of $\mathrm{H}^n(F_x)_\lambda$ is supported in at most one Hodge bidegree. If all singularities of $X_0$ are isolated, we say that $X_0$ has concentrated mixed Hodge spectrum when the sum of their spectra has this property.
\end{defn}

The mixed spectrum can be calculated explicitly for isolated quasihomogeneous singularities via the combinatorics of the weights, see \cite{steenbrink1977intersection,steenbrink1985mixed,KerrLazaToAppear}.
\begin{ex}
	\label{example: A-D-E}
For ordinary double points $(X,x)$ of dimension $n$, the mixed Hodge spectrum is given by
\[\sigma(X, x) = \begin{cases}
[(\frac{n}{2}+\frac{1}{2}, n)], & n \text{ even}\\
[(\frac{n+1}{2}, n+1)], & n \text{ odd}
\end{cases}\]	
For cusp singularity $(Y,y)$ of dimension $n$, i.e. locally given by $y_0^2+y_1^2+\cdots+y_n^3=0$, then the mixed Hodge spectrum is given by
\[\sigma(Y, y) = \begin{cases}
[(\frac{n}{2}+\frac{1}{3}, n)] + [(\frac{n}{2}+\frac{2}{3}, n)], & n \text{ even}\\
[(\frac{n-1}{2}+\frac{5}{6}, n)] + [(\frac{n+1}{2}+\frac{1}{6}, n)], & n \text{ odd}
\end{cases}\]

Thus any collection of ordinary double points and cusp singularities in the same dimension has concentrated summed mixed Hodge spectrum: for each fixed eigenvalue and weight, the displayed spectra have only one Hodge index. Consequently, the hypothesis of Theorem~\ref{main-thm} is satisfied if all singularities of the central fiber are of these types.
\end{ex}

\subsection{Deformation of the top form under a smoothing}
Let $\Omega$ be the nowhere vanishing relative holomorphic top form in
Assumption~\ref{assum: main}, and denote by $[\Omega_t]$ the cohomology
class of $\Omega_t=\Omega|_{X_t}$ in $\mathrm{H}^n(X_t,\CC)$.
We first show that its restriction to the Milnor fiber at every isolated
singular point is nonzero.

\begin{prop}\label{non-triv-T}
In the setting of Assumption~\ref{assum: main}, let $p\in X_0$ be an
isolated singular point. For $0<|t|\ll1$, let
\[
\rho_{p,t}:\mathrm{H}^n(X_t,\CC)\longrightarrow
\mathrm{H}^n(F_{p,t},\CC)
\]
be the restriction to the Milnor fiber. Then
\[
\rho_{p,t}([\Omega_t])\neq0.
\]
In particular, if $\mathbb V_0$ is a polarized $\CC$-VHS summand
containing the Hodge line $F^n$, then the image of
$\Hlim(\mathbb V_0)$ under the local vanishing-cycle map $\nu_p$
is nonzero.
\end{prop}

\section{Monodromy action and the period integral}\label{sec-3}
The goal of this section is to prove Proposition~\ref{non-triv-T}. We work over the disk in Assumption~\ref{assum: main}, and write it simply as $f:\calX\to\Delta$. Thus $\calX$ is smooth and $\omega_{\calX/\Delta}$ is generated by the relative volume form $\Omega$.

\subsection{A local residue proof of Proposition~\ref{non-triv-T}}
The following elementary point is where the Calabi--Yau condition enters the
argument. It is not enough to take an arbitrary section of the relative
dualizing sheaf: the section must be a local generator at the chosen isolated
singular point. This fails, for instance, for a degeneration of genus two
curves with a separating node, where the limiting holomorphic one-forms vanish
in the local dualizing direction at the node.

We record the local input for isolated hypersurface singularities. Let
$g:(Y,p)\to (\Delta,0)$ be a holomorphic function germ on a smooth
$(n+1)$-fold, and assume that $p$ is an isolated critical point of $g$.
For $0<|t|\ll 1$, denote by $F_t=g^{-1}(t)\cap B_{\epsilon}(p)$ the Milnor
fiber.

\begin{lem}\label{lem:local-residue-nonzero}
Let $\omega$ be a local generator of the relative dualizing sheaf
$\omega_{Y/\Delta}$ near $p$. Then the de Rham class of
$\omega|_{F_t}$ in $\mathrm{H}^n(F_t,\CC)$ is nonzero for $0<|t|\ll 1$.
\end{lem}

\begin{proof}
Choose local coordinates $z_0,\ldots,z_n$ on $Y$ centered at $p$, and
write $\eta=dz_0\wedge\cdots\wedge dz_n$. Since $\omega$ generates
$\omega_{Y/\Delta}$, its restriction to a smooth fiber is
\[
\omega|_{F_t}=\operatorname{Res}_{g=t}\frac{u\eta}{g-t}
\]
for a unit $u\in\mathcal O_{Y,p}^{\times}$.

The Brieskorn lattice
\[
H''_g:=\frac{\Omega^{n+1}_{Y,p}}
{dg\wedge d\Omega^{n-1}_{Y,p}}
\]
is a free $\CC\{t\}$-module of rank equal to the Milnor number, with
$t[\alpha]=[g\alpha]$ \cite{Bri70,Seb70}. Its fiber over $0$ is
\[
\frac{H''_g}{tH''_g}=
\frac{\Omega^{n+1}_{Y,p}}
{g\Omega^{n+1}_{Y,p}+dg\wedge d\Omega^{n-1}_{Y,p}}.
\]
We claim that the image of $u\eta$ in this quotient is nonzero.
Otherwise, there would exist $v\in\mathcal O_{Y,p}$ and
$\beta\in\Omega^{n-1}_{Y,p}$ such that
\[
u\eta=gv\eta+dg\wedge d\beta.
\]
Writing the coefficient of $dg\wedge d\beta$ relative to $\eta$ shows
that
\[
u\in(g,J_g)\subset\mathfrak m_p,
\qquad
J_g=(\partial g/\partial z_0,\ldots,\partial g/\partial z_n),
\]
contrary to the fact that $u$ is a unit.

Choose a $\CC\{t\}$-basis $e_1,\ldots,e_\mu$ of $H''_g$, and write
$[u\eta]=\sum_i a_i(t)e_i$. The nonvanishing modulo $t$ means that
$a_i(0)\neq0$ for some $i$. After shrinking the disk, this coefficient
remains nonzero. The Brieskorn--Gauss--Manin comparison identifies the
specialization at $t\neq0$ with the de Rham class of the displayed
residue form on $F_t$ \cite{Bri70}. Hence
$[\omega|_{F_t}]\neq0$ for $0<|t|\ll1$.
\end{proof}

\begin{rmk}
For a general isolated hypersurface singularity, one must not identify
$H''_g/tH''_g$ with the Jacobian algebra
$(\mathcal O_{Y,p}/J_g)\eta$. The argument above uses the actual
$t$-fiber and does not require quasihomogeneity.
\end{rmk}

\begin{proof}[Proof of Proposition~\ref{non-triv-T}]
Since $\calX$ is smooth and $X_0=f^{-1}(0)$ is a Cartier divisor, an
isolated singular point $p\in X_0$ is an isolated critical point of the
local function $f$ on $\calX$. By Assumption~\ref{assum: main},
$\Omega$ generates $\omega_{\calX/\Delta}$ at $p$.
Lemma~\ref{lem:local-residue-nonzero}, applied to the germ
\[
f:(\calX,p)\longrightarrow(\Delta,0),
\]
therefore gives
\[
\rho_{p,t}([\Omega_t])=[\Omega_t|_{F_{p,t}}]\neq0.
\]

The restriction maps $\rho_{p,t}$ form a morphism of local systems on
the punctured disk. Since $[\Omega_t]$ belongs to $\mathbb V_{0,t}$,
this morphism is nonzero on $\mathbb V_0$. Passing to the limiting
mixed Hodge structures gives $\nu_p|_{\Hlim(\mathbb V_0)}$, which has
the same underlying linear map after choosing a lift to the universal
cover of the punctured disk. Its image is consequently nonzero.
\end{proof}

\begin{rmk}\label{rmk:nearby-limit-nonvanishing}
Proposition~\ref{non-triv-T} concerns the top form on nearby fibers.
It implies nonvanishing of the limiting vanishing-cycle map on
$\Hlim(\mathbb V_0)$, but does not assert that a generator of the
limiting top Hodge line has nonzero image. Indeed, for an ordinary
double point on a Calabi--Yau threefold, the local vanishing cohomology
has type $(2,2)$, so the map on $F^3$ of the limiting mixed Hodge
structure is zero. Only the nonzero image of the entire summand is
needed below.
\end{rmk}

\begin{rmk}\label{rmk:separating-node}
The preceding proof would be false without the assumption that the relative
volume form generates $\omega_{\calX/\Delta}$ at the chosen singular point. For
a genus-two curve degeneration with a separating node, the local model is $xy=t$, but
a limiting holomorphic one-form can have local expression
$x\,d x/x=d x$ or $y\,d x/x$ near the node. Its class on the annular Milnor
fiber is zero, and the corresponding Picard--Lefschetz transformation is
trivial on $\mathrm H^1$ because the vanishing loop is zero in the global homology. The relative
Calabi--Yau condition in Assumption~\ref{assum: main} rules out this phenomenon
by forcing the relative volume form to be a local generator of the dualizing
sheaf at $p$.
\end{rmk}

\subsection{Period integral approach to ODP}
When all singularities of $X_0$ are ordinary double points, nonvanishing can also be seen directly from the period integrals. This is related to the discussion in \cite{rollenske2009smoothing}; see also \cite[Corollary~3.5]{wang1997incompleteness} for the Clemens--Schmid approach. We give the local residue calculation below and recall its differential geometric interpretation in dimension three.

\begin{prop}
	\label{prop: ordinary double points}
Let $\calX$ be a family of Calabi--Yau varieties satisfying Assumption~\ref{assum: main}. Assume, in addition, that all singularities of $X_0$ are ordinary double points. Then 
\begin{enumerate}
	\item When $n$ is even, then $N=0$, $T^{ss}\neq Id$ and the $(-1)$-eigenspace $(\Hlim^n(X_t))^{T=-Id}\neq 0$ has pure Hodge structure with Hodge numbers $h^{p,q}\neq 0$ if and only if $(p,q)=(\frac{n}{2},\frac{n}{2})$.
	\item When $n$ is odd, then $N\neq 0$, $T^{ss}=Id$ and $\Ima(N)=W_{n-1}(\Hlim^n(X_t))\neq 0$ has pure Hodge structure with Hodge numbers $h^{p,q}\neq 0$ if and only if $(p,q)=(\frac{n-1}{2},\frac{n-1}{2})$.
\end{enumerate}
\end{prop}

	\begin{proof}[Proof of Proposition~\ref{prop: ordinary double points}]
Put $H=\Hlim^n(X_t)$ and
$H_{\rm van}=\bigoplus_i\mathrm H^n(F_{p_i})$.
The vanishing-cycle sequence and the invariant-cycle theorem give
\[
\ker(H\longrightarrow H_{\rm van})=H^T.
\]
By Proposition~\ref{non-triv-T}, the image of this map is nonzero.
For an ordinary double point, $\mathrm H^n(F_{p_i})$ is one-dimensional,
of type $(n/2,n/2)$ with monodromy $-1$ if $n$ is even, and of type
$((n+1)/2,(n+1)/2)$ with monodromy $1$ if $n$ is odd.

When $n$ is even, the quotient $H/H^T$ has monodromy $-1$, while $T$
acts as the identity on $H^T$. Hence $T$ is semisimple and its nonzero
$(-1)$-eigenspace embeds into $H_{\rm van}$. The assertion follows.

When $n$ is odd, $(T-\id)H\subset H^T$, so $(T-\id)^2=0$ and
$T^{ss}=\id$. Since $H/H^T$ is nonzero, $N\neq0$. The monodromy weight
filtration is then
\[
0\subset W_{n-1}=\Ima N\subset W_n=\ker N\subset W_{n+1}=H.
\]
The quotient $\Gr^W_{n+1}H=H/\ker N$ embeds into $H_{\rm van}$,
and is therefore of type $((n+1)/2,(n+1)/2)$. The isomorphism of pure
Hodge structures
\[
N:\Gr^W_{n+1}H\xrightarrow{\sim}(\Gr^W_{n-1}H)(-1)
\]
shows that $\Ima N$ has type $((n-1)/2,(n-1)/2)$.
\end{proof}
 For $n\geq2$, we now look at corresponding mixed Hodge structures on $\mathrm{H}^n(X_0)$. Let $\widetilde{X_0}\to X_0$ be the resolution of singularity of $X_0$ such that the exceptional divisors over $p_i$ are all $n-1$ dimensional smooth quadrics $E_i$. Then we have the following exact sequence of mixed Hodge structures
		\begin{equation*}
			\cdots \to \bigoplus_{i=1}^k \mathrm{H}^{n-1}(E_i)\to \mathrm{H}^n(X_0)\to \mathrm{H}^n(\widetilde{X_0})\to \bigoplus_{i=1}^k \mathrm{H}^{n}(E_i)\to \mathrm{H}^{n+1}(X_0)\to \cdots
		\end{equation*}
	From the cohomology of quadrics, we know that $\mathrm{H}^l(E_i)\neq 0$ if and only if $l$ is even and $0\leq l\leq 2n-2$. In this case, the Hodge numbers are concentrated in the middle degree. When $n$ is even, then $\mathrm{H}^n(X_0)$ is pure of weight $n$ because $\mathrm{H}^{n-1}(E_i)=0$. So this also implies $\Hlim^n(X_t)$ is pure of weight $n$ and $N=0$.

	When $n\geq3$ is odd, the same exact sequence gives
\[
W_{n-1}\mathrm H^n(X_0)\simeq
\operatorname{coker}\left(\mathrm H^{n-1}(\widetilde X_0)
\longrightarrow\bigoplus_i\mathrm H^{n-1}(E_i)\right).
\]
Each $\mathrm H^{n-1}(E_i)$ has dimension two. The normal bundle of
$E_i$ is $\calO_{E_i}(-1)$, where $\calO_{E_i}(1)$ is the hyperplane
bundle on the quadric. Thus the restrictions of
$c_1(\calO_{\widetilde X_0}(E_i))^{(n-1)/2}$ span one independent
class on each $E_i$. By the invariant-cycle theorem and
Proposition~\ref{prop: ordinary double points}, the dimension $a$ of
the displayed cokernel satisfies
\[
0<a=\rank N\leq k.
\]
Moreover, $\mathrm H^{n+1}(X_0)$ embeds into
$\mathrm H^{n+1}(\widetilde X_0)$ and is pure of weight $n+1$.
The vanishing-cycle sequence also gives a relation among the images
of the local vanishing classes: if
$\rho_t([\Omega_t])=(c_i(t)e_i)_i$, where $e_i$ generates
$\mathrm H^n(F_{p_i,t})$, then every $c_i(t)$ is nonzero and
\[
\sum_i c_i(t)\,\partial(e_i)=0
\quad\hbox{in }\mathrm H^{n+1}(X_0).
\]
Here $\partial$ is the connecting map in
\eqref{equation: vanishing sequences}. These relations are related
to the smoothing conditions studied by Friedman
\cite{friedman1986simultaneous} in dimension three and by
Rollenske--Thomas \cite{rollenske2009smoothing} in higher odd
dimensions.

\begin{lem}
		\label{lemma: vanishing cycles in ODP}
		In the setting of Proposition \ref{prop: ordinary double points} and let $\Omega_t$ be the generator of $\mathrm{H}^0(X_t, K_{X_t})$ and $\delta_i$ be the vanishing cycle near $p_i$, then $\int_{\delta_i} \Omega_t \neq 0$ for $0<|t|<\epsilon$ with $\epsilon$ small enough. In fact $\int_{\delta_i} \Omega_t \approx c \cdot t^{\frac{n-1}{2}}$ for some constant $c \neq 0$ when $\Omega_t$ is the generator of $f_*K_{\calX/\Delta}$.
	\end{lem}

	\begin{proof}[Proof of Lemma \ref{lemma: vanishing cycles in ODP}]
Choose local coordinates in which $f=x_0^2+\cdots+x_n^2=t$.
Since $\Omega$ generates the relative canonical bundle, write
\[
\Omega_t=\operatorname{Res}
\frac{u(x)\,dx_0\wedge\cdots\wedge dx_n}
{x_0^2+\cdots+x_n^2-t},\qquad u(0)\ne0.
\]
On a sector with a choice of $\sqrt t$, take the vanishing sphere
$\delta_t=\sqrt t\,S^n$. The substitution $x=\sqrt t\,y$ gives
\[
\int_{\delta_t}\Omega_t
=t^{\frac{n-1}{2}}\left(c\,u(0)+O(|t|^{1/2})\right),
\qquad
c=\int_{S^n}\operatorname{Res}
\frac{dy_0\wedge\cdots\wedge dy_n}{\sum y_i^2-1}\ne0.
\]
Indeed, up to orientation, the model residue restricts to one half of
the standard volume form on $S^n$. This proves the asserted asymptotic
and nonvanishing in every dimension.

For nodal Calabi--Yau threefolds, the special Lagrangian construction and
the estimates for the holomorphic volume form in
\cite[Corollary~A.2]{HeinSun2017} and
\cite[Theorem~1.1 and Lemma~4.3]{collins2023special} give a
differential geometric interpretation of this nonzero period.
\end{proof}

\section{Tensor-product decomposition of limiting MHS}\label{sec-4}
To prove Theorem~\ref{main-thm}, we argue by contradiction. Suppose that $f$ is {\it non-rigid}. Thus it can be deformed to a family of polarized Calabi-Yau manifolds over a product base $S \times B$, where $B$ is some quasiprojective curve. The enlarged family is non-isotrivial in the $B$ direction for general $s\in S$.

Denote by the enlarged family as $g\colon (\calZ, \calL_Z)\to S \times B$. Let $\mathbb{V}$ be the primitive part of the $\mathbb{C}$-VHS $R^ng_*\mathbb{C}_{\calZ}$. Then $\mathbb{V}$ is a polarized $\mathbb{C}$-VHS of Calabi-Yau type over $S \times B$.
Let $\mathbb{V}_0$ be the irreducible factor of $\mathbb{V}$ containing the deepest rank one Hodge filtration $F^nV$ of $\mathbb{V}$. Thus $\mathbb{V}_0$ is again a $\mathbb{C}$-VHS of Calabi-Yau type over $S \times B$.

\begin{prop}[Tensor-product decomposition]\label{Deligne}
There exist polarized $\mathbb{C}$-VHS $\mathbb{V}_1$ over $S$ and $\mathbb{V}_2$ over $B$, and a Hodge isometry 
\begin{equation}
  \label{tensor-decomp}
\mathbb{V}_0 \simeq p^*_1\mathbb{V}_1 \otimes p^*_2\mathbb{V}_2.
\end{equation}
\end{prop}

\begin{rmk}
Proposition~\ref{Deligne} is essentially due to Deligne \cite{deligne1987theoreme}, a detailed proof can be found in \cite[Proposition 3.3]{viehweg2005complex}. Note that the quasi-projectivity of $S$ and $B$ are crucial here since the semisimplicity of $\mathbb{V}_i$ are used in the proof.
\end{rmk}

Therefore, at the boundary point $0\in\overline S\setminus S$ and the original base point $b_0\in B$, one has the isomorphism
\begin{align}
  \label{H_lim}
\mathrm{H}_{\mathrm{lim}}(\mathbb{V}_{0,(0,b_0)}) \simeq \mathrm{H}_{\mathrm{lim}}(\mathbb{V}_{1,0}) \otimes \mathrm{H}(\mathbb{V}_{2,b_0})
\end{align}
between mixed Hodge structures.

\begin{proof}[Proof of Theorem~\ref{main-thm}]
Suppose that the family is non-rigid, and write the decomposition \eqref{H_lim} as
\[
H=A\otimes K,\qquad
A=\Hlim(\mathbb{V}_{1,0}),\quad K=\mathrm{H}(\mathbb{V}_{2,b_0}).
\]
The monodromy acts as $T=T_A\otimes\Id_K$. By \eqref{eq:total-kernel}, the total vanishing-cycle image of $H$ is
\begin{equation}\label{eq:total-product-image}
\Ima(\nu_\Sigma|_H)\simeq H/H^T
\simeq (A/A^{T_A})\otimes K
\end{equation}
as mixed Hodge structures with semisimple monodromy. By the nearby-cycle comparison in Proposition~\ref{non-triv-T}, this image is nonzero, since $H$ is the limiting subvariation containing the top Hodge line. Hence $A/A^{T_A}\ne0$.

The non-isotriviality in the $B$ direction and the local Torelli property for Calabi--Yau manifolds imply that the Higgs field of $\mathbb{V}_2$ is nonzero on its highest Hodge line at a general point of $B$. Thus $K$ has two nonzero adjacent Hodge components, say $K^{c,d}$ and $K^{c-1,d+1}$; their dimensions are constant on $B$. Choose a nonzero Hodge--Deligne component of $A/A^{T_A}$ of type $(r,s)$ and monodromy eigenvalue $e^{2\pi\sqrt{-1}\lambda}$. Formula~\eqref{eq:total-product-image} then gives two nonzero components of types
\[
(r+c,s+d),\qquad (r+c-1,s+d+1)
\]
with the same eigenvalue and the same weight. Since $\Ima(\nu_\Sigma|_H)$ is a mixed Hodge substructure of $H_\Sigma$, this contradicts concentration of $\sigma(X_0)$.
\end{proof}

\begin{rmk}
The use of the total vanishing-cycle map is essential here. Nonvanishing at a single isolated point does not imply that the image of its local map retains the whole factor $K$. Identity~\eqref{eq:total-kernel} supplies this assertion for the total map. Accordingly, the hypothesis concerns the sum of all local spectra, not merely concentration at each point separately.
\end{rmk}

\begin{ex}[A nonrigid family with an isolated singularity]\label{example:nonrigid-cone}
Let $B\subset\mathbb A^1$ be the open set of parameters $b$ for which
the binary quintic $u^5+v^5+bu^3v^2$ is square-free, and put
\[
S=\{s\in\mathbb A^1\setminus\{0\}\mid s^2\in B\}.
\]
The equation
\[
x_0^5+x_1^5+x_2^5+x_3^5+s^2x_0^3x_1^2+bx_2^3x_3^2+sw^5=0
\quad\text{in }\mathbb P^4
\]
defines a smooth family of quintic threefolds over $S\times B$.
The restriction to $b=0$ is non-isotrivial and admits a nontrivial
deformation in the $b$-direction, so it is nonrigid.

At $s=0$ and $b=0$, the central fiber is the projective cone over the
smooth Fermat quintic surface
\[
\{x_0^5+x_1^5+x_2^5+x_3^5=0\}\subset\mathbb P^3,
\]
with a unique singularity at the vertex $p=[0:0:0:0:1]$.
The total space of the $b=0$ family is smooth near the central fiber,
and its relative canonical bundle is trivial.
The local mixed Hodge spectrum at $p$ contains
\[
[(16/5,3)]\quad\text{and}\quad[(11/5,3)]
\]
by the quasihomogeneous formula
\cite{steenbrink1977intersection,KerrLazaToAppear}.
These have the same monodromy eigenvalue $e^{2\pi\sqrt{-1}/5}$ and
weight $3$, but distinct Hodge indices $3$ and $2$.
Thus the spectrum is not concentrated. This example shows that
isolated singularities alone do not imply rigidity, even for a
non-isotrivial family.
\end{ex}

\section{The infinitesimal version}\label{sec-5}
We now give an infinitesimal version of Theorem~\ref{main-thm} in
terms of flat endomorphisms. Let $\mathbb W$ be the primitive part of
$R^nf_*\CC$, and let $\mathbb W_0$ be the isotypic component containing
the top Hodge line $F^n\mathbb W$. Put
\[
E_0=\operatorname{End}_{\pi_1(S)}(\mathbb W_0).
\]
This is the algebra of global flat endomorphisms of $\mathbb W_0$,
endowed with its natural weight-zero Hodge structure.

\begin{thm}\label{thm:infinitesimal}
Under the assumptions of Theorem~\ref{main-thm}, one has
\[
E_0^{-1,1}=0.
\]
Consequently, every
$A\in\operatorname{End}_{\pi_1(S)}(\mathbb W)^{-1,1}$ annihilates
the top Hodge line $F^n\mathbb W$.
\end{thm}

\begin{proof}
By semisimplicity and the theorem of the fixed part, $\mathbb W_0$
is a subvariation and
\[
E_0\simeq\operatorname{End}_{\CC}(M),
\]
where $M$ is the multiplicity space of the corresponding irreducible
local system. In particular, $E_0$ is a simple algebra. Indeed, the
Hodge action on the full flat endomorphism algebra fixes its primitive
central idempotents, since these form a finite set. The evaluation map
$E_0\otimes\mathbb W_0\longrightarrow\mathbb W_0$ is a morphism of
variations of Hodge structures.

Write $\operatorname{Sing}(X_0)=\{p_1,\ldots,p_r\}$ and put
\[
H_\Sigma=\bigoplus_{i=1}^r\mathrm H^n(F_{p_i},\CC),\qquad
H=\Hlim(\mathbb W_0),\qquad
I_0=\operatorname{Im}\bigl(H\xrightarrow{\nu_\Sigma}H_\Sigma\bigr).
\]
The total vanishing-cycle sequence and the invariant-cycle theorem
give, as in \eqref{eq:total-kernel},
\[
I_0\simeq H/H^T.
\]
By Proposition~\ref{non-triv-T}, $I_0\ne0$. Every element of $E_0$
commutes with $T$ and therefore acts on this quotient. The kernel of
\[
E_0\longrightarrow\operatorname{End}_{\CC}(I_0)
\]
is a two-sided ideal of the simple algebra $E_0$. Since the identity
acts nontrivially on $I_0$, this action is faithful.

Let $A\in E_0^{-1,1}$. Passing the evaluation map to limiting mixed
Hodge structures and then to the quotient, functoriality of the
Deligne splitting gives
\[
A:I_{0,\lambda}^{p,q}\longrightarrow I_{0,\lambda}^{p-1,q+1}.
\]
Here $I_{0,\lambda}^{p,q}$ denotes the Hodge--Deligne component with
semisimple monodromy eigenvalue $e^{2\pi\sqrt{-1}\lambda}$.
Thus $A$ preserves the eigenvalue and the weight $p+q$.
Since $I_0$ is a mixed Hodge substructure of $H_\Sigma$, strictness
implies that its mixed Hodge spectrum is concentrated. The source
and target of the displayed map cannot both be nonzero. Hence $A$
acts trivially on $I_0$, and faithfulness gives $A=0$.

Finally, every flat endomorphism of $\mathbb W$ preserves its
isotypic decomposition. An endomorphism of type $(-1,1)$ therefore
vanishes on $\mathbb W_0$, and in particular on $F^n\mathbb W$.
\end{proof}

The infinitesimal statement also gives another proof of our rigidity
criterion.

\begin{cor}\label{cor:infinitesimal-rigidity}
Under the assumptions of Theorem~\ref{main-thm}, the family
$f\colon\calX\to S$ is rigid.
\end{cor}

\begin{proof}
It suffices to consider a deformation $\widetilde f$ of $f$ over
$S\times B$, where $(B,b_0)$ is a smooth pointed curve. Locally on
$B$, the isotypic component $\mathbb W_0$ extends to a subvariation
$\widetilde{\mathbb W}_0$ containing the top Hodge line. By the
theorem of the fixed part, the spaces
\[
E_{0,b}=\operatorname{End}_{\pi_1(S)}
\bigl(\widetilde{\mathbb W}_0|_{S\times\{b\}}\bigr)
\]
form a weight-zero variation on $B$: they are the Hodge structures
on the flat endomorphisms invariant under $\pi_1(S)$, or equivalently
the endomorphism Hodge structures of the multiplicity spaces
described in \cite[\S4]{peters90}. Their Hodge numbers are locally
constant. Thus Theorem~\ref{thm:infinitesimal} implies
$E_{0,b}^{-1,1}=0$ for every $b$ near $b_0$.

By \cite[Theorem~3.2(i)]{peters90}, each infinitesimal deformation
in the $B$-direction induces a global flat endomorphism of type
$(-1,1)$ of the primitive variation on $S\times\{b\}$. It preserves
the isotypic decomposition, so its restriction to
$\widetilde{\mathbb W}_0|_{S\times\{b\}}$ is zero. In particular, it
annihilates the top Hodge line. Infinitesimal Torelli for
Calabi--Yau manifolds then implies that every Kodaira--Spencer map
in the $B$-direction is zero. Hence the classifying map is locally
constant in that direction, and $f$ is rigid.
\end{proof}

\section{Relation to C.-L. Wang's result on Weil--Petersson metric}\label{sec-6}
In this section, we discuss the relation between our main result Theorem~\ref{main-thm} and C.-L. Wang's result \cite[Proposition~2.3]{wang1997incompleteness} on incompleteness of Weil--Petersson metric for Calabi--Yau families. We specialize to the case of Calabi--Yau threefold with ordinary double points as singularities. Then C.-L. Wang's \cite[Proposition~2.3]{wang1997incompleteness} implies that the Weil--Petersson metric is incomplete along any smoothing for the singular fiber $X_0$.

On the other hand, if the family $f\colon (\calX, \calL)\to S$ is non-rigid, then by Proposition~\ref{Deligne}, one has the tensor-product decomposition of $\mathbb{C}$-VHS
\[\mathbb{V}_0 \simeq p^*_1\mathbb{V}_1 \otimes p^*_2\mathbb{V}_2.\]
Here $\mathbb{V}_0$ is the irreducible factor of the primitive part of $R^3f_*\mathbb{C}_{\calX}$ containing the deepest Hodge filtration. By the classfication of all possible nontrivial tensor product decomposition of weight-3 complex Hodge structures of Calabi--Yau type in \cite[Proposition~14]{viehweg2005geometry}, the possible period domains corresponding to $\CC$-VHS on $\mathbb{V}_1$ and $\mathbb{V}_2$ are either the complex ball $\mathrm{SU}(1,n)/\mathrm{S}(\mathrm{U}(1)\times \mathrm{U}(n))$ or the type IV symmetric domain $\mathrm{SO}(2,n)/\mathrm{SO}(2)\times \mathrm{SO}(n)$. In both cases, the Hodge metric on the period domains are known to be complete for boundary points with nontrivial unipotent monodromy.

By Proposition~\ref{prop: ordinary double points} or \cite[Theorem~1.5(1)]{wang2003quasi}, the monodromy operator $T$ on $\mathrm{H}_{\mathrm{lim}}^3(X_t)$ satisfies that $N \neq 0$ and $T^{ss}=Id$. So the boundary point corresponding to the singular fiber $X_0$ is at infinite distance with respect to the Hodge metric on the period domain of $\mathbb{V}_1$. Since the two period domains corresponding to $\mathbb{V}_1$ and $\mathbb{V}_2$ are both Hermitian symmetric domainis, the Hodge metric is the same as the Weil--Petersson metric up to a constant multiple. Therefore, the Weil--Petersson metric on $S$ is complete along the smoothing for the singular fiber $X_0$. This contradicts to C.-L. Wang's result \cite[Proposition~2.3]{wang1997incompleteness}.

Notice that the metric completeness argument here works only for Calabi--Yau threefolds with ordinary double points as singularities. For more general singularities or higher dimensional Calabi--Yau varieties, the monodromy operator may not have nontrivial unitpotent part, or the decomposition factors of the tensor-product decomposition may not correspond to Hermitian symmetric domains. So the metric completeness argument may not introduce the desired contradiction. But more detailed study of the curvature behavior of Weil--Petersson metric near singular fibers may lead to alternative proof of Theorem~\ref{main-thm} other than the mixed Hodge structure approach.


\begin{thebibliography}{CGPY23}

\bibitem[Ara71]{arakelov1971}
S.~J. Arakelov, \emph{Families of algebraic curves with fixed degeneracies},
  Izvestiya Akademii Nauk SSSR. Seriya Matematicheskaya \textbf{35} (1971),
  1269--1293.

\bibitem[Bri70]{Bri70}
E.~Brieskorn, \emph{Die Monodromie der isolierten Singularit\"aten von
Hyperfl\"achen}, Manuscripta Math. \textbf{2} (1970), 103--161.

\bibitem[CGPY23]{collins2023special}
T.~C. Collins, S.~Gukov, S.~Picard, and S.-T. Yau, \emph{Special {Lagrangian}
  cycles and {Calabi}--{Yau} transitions}, Communications in Mathematical
  Physics \textbf{401} (2023), no.~1, 769--802.

\bibitem[CHSZ24]{CHSZ}
K.~Chen, T.~Hu, R.~Sun, and K.~Zuo, \emph{On the distribution of non-rigid
  families in the moduli spaces}, 2024, arXiv:2408.11604.

\bibitem[Del87]{deligne1987theoreme}
P.~Deligne, \emph{Un th{\'e}oreme de finitude pour la monodromie}, Discrete
  groups in geometry and analysis: papers in honor of GD Mostow on his sixtieth
  birthday, Springer, 1987, pp.~1--19.

\bibitem[Fal83]{faltings83}
G.~Faltings, \emph{Arakelov's theorem for abelian varieties}, Invent. Math.
  \textbf{73} (1983), no.~3, 337--347. \MR{718934}

\bibitem[FL24]{friedman2024deformations}
R.~Friedman and R.~Laza, \emph{Deformations of {Calabi}--{Yau} varieties with
  k-liminal singularities}, Forum of Mathematics, Sigma, vol.~12, Cambridge
  University Press, 2024, p.~e59.

\bibitem[FL25]{friedman2025deformations}
\bysame, \emph{Deformations of singular {Fano} and {Calabi--Yau} varieties},
  Journal of Differential Geometry \textbf{131} (2025), no.~1, 65--131.

\bibitem[Fri86]{friedman1986simultaneous}
R.~Friedman, \emph{Simultaneous resolution of threefold double points},
  Mathematische Annalen \textbf{274} (1986), no.~4, 671--689.

\bibitem[HS17]{HeinSun2017}
H.~Hein and S.~Sun, \emph{{Calabi}--{Yau} manifolds with isolated conical
  singularities}, Publications Math\'ematiques de l'IH\'ES \textbf{126} (2017),
  73--130.

\bibitem[KL]{KerrLazaToAppear}
M.~Kerr and R.~Laza, \emph{Hodge theory of degenerations, {II}: vanishing
  cohomology and geometric applications}, IMSA Hodge theory proceedings, to
  appear.

\bibitem[KL21]{KerrLaza2021}
\bysame, \emph{Hodge theory of degenerations, (i): consequences of the
  decomposition theorem}, Selecta Mathematica (New Series) \textbf{27} (2021),
  no.~4, Paper No. 71, 48, With an appendix by M. Saito.

\bibitem[LTYZ05]{LTYZ05}
K.~Liu, A.~Todorov, S.-T. Yau, and K.~Zuo, \emph{Shafarevich's conjecture for
  {CY} manifolds. {I}}, Q. J. Pure Appl. Math. \textbf{1} (2005), no.~1,
  28--67. \MR{2155142}

\bibitem[Mil68]{milnor1968singular}
J.~W. Milnor, \emph{Singular points of complex hypersurfaces}, no.~61,
  Princeton University Press, 1968.

\bibitem[Par68]{parshin1968}
A.~N. Parshin, \emph{Algebraic curves over function fields. {I}}, Izvestiya
  Akademii Nauk SSSR. Seriya Matematicheskaya \textbf{32} (1968), 1191--1219.

\bibitem[Pet90]{peters90}
C.~A.~M. Peters, \emph{Rigidity for variations of {H}odge structure and
  {A}rakelov-type finiteness theorems}, Compositio Math. \textbf{75} (1990),
  no.~1, 113--126. \MR{1059957}

\bibitem[PS08]{peters2008mixed}
C.~Peters and J.~Steenbrink, \emph{Mixed {Hodge} structures}, Springer, 2008.

\bibitem[RT09]{rollenske2009smoothing}
S.~Rollenske and R.~Thomas, \emph{Smoothing nodal {Calabi}--{Yau} n-folds},
  Journal of Topology \textbf{2} (2009), no.~2, 405--421.

\bibitem[Sai93]{saito-ab}
M.-H. Sait{o}, \emph{Classification of nonrigid families of abelian varieties},
  Tohoku Math. J. (2) \textbf{45} (1993), no.~2, 159--189. \MR{1215923}

\bibitem[Seb70]{Seb70}
M.~Sebastiani, \emph{Preuve d'une conjecture de Brieskorn},
Manuscripta Math. \textbf{2} (1970), 301--308.

\bibitem[Ste77]{steenbrink1977intersection}
J.~Steenbrink, \emph{Intersection form for quasi-homogeneous singularities},
  Compositio Mathematica \textbf{34} (1977), no.~2, 211--223.

\bibitem[Ste85]{steenbrink1985mixed}
\bysame, \emph{Mixed {H}odge structures and vanishing cohomology}, Compositio
  Mathematica \textbf{57} (1985), 321--344.

\bibitem[SYZ]{sun2026nonrigidCY}
R.~Sun, C.~Yu, and K.~Zuo, \emph{Nonrigid family of {Calabi}--{Yau} manifolds
  and their moduli}.

\bibitem[SZ91]{saito-zucker}
M.-H. Sait{o} and S.~Zucker, \emph{Classification of nonrigid families of
  {$K3$} surfaces and a finiteness theorem of {A}rakelov type}, Math. Ann.
  \textbf{289} (1991), no.~1, 1--31. \MR{1087233}

\bibitem[VZ04]{VZ04Ara}
E~Viehweg and K~Zuo, \emph{A characterization of certain {S}himura curves in
  the moduli stack of abelian varieties}, J. Differential Geom. \textbf{66}
  (2004), no.~2, 233--287. \MR{2106125}

\bibitem[VZ05]{viehweg2005complex}
E.~Viehweg and K.~Zuo, \emph{Complex multiplication, {Griffiths}--{Yukawa}
  couplings, and rigidity for families of hypersurfaces}, Journal of Algebraic
  Geometry \textbf{14} (2005), no.~3, 481--528.

\bibitem[VZ06]{viehweg2005geometry}
\bysame, \emph{Geometry and arithmetic of non-rigid families of {C}alabi-{Y}au
  3-folds; questions and examples}, Mirror symmetry. {V}, AMS/IP Stud. Adv.
  Math., vol.~38, Amer. Math. Soc., Providence, RI, 2006, pp.~351--360.
  \MR{2282966}

\bibitem[Wan97]{wang1997incompleteness}
C.-L. Wang, \emph{On the incompleteness of the {Weil--Petersson} metric along
  degenerations of {Calabi--Yau} manifolds}, Mathematical Research Letters
  \textbf{4} (1997), no.~1, 157--171.

\bibitem[Wan03]{wang2003quasi}
\bysame, \emph{{Quasi-Hodge} metrics and canonical singularities}, Mathematical
  Research Letters \textbf{10} (2003), no.~1, 57--70.

\bibitem[Zha04]{zhang04}
Y.~Zhang, \emph{Rigidity for families of polarized {C}alabi-{Y}au varieties},
  J. Differential Geom. \textbf{68} (2004), no.~2, 185--222. \MR{2144247}

\end{thebibliography}

\end{document}